\documentclass[xypic,amscd,syntonly,amssymb]{amsart}

\usepackage{graphicx}
\usepackage[all]{xy}
\usepackage{amssymb}
\usepackage{amsmath}

\usepackage{epsfig}
\theoremstyle{plain}
\newtheorem{Thm}{Theorem}
\newtheorem{Prop}[Thm]{Proposition}
\newtheorem{Cor}[Thm]{Corollary}

 \theoremstyle{definition}
\theoremstyle{remark}

\numberwithin{equation}{section}

\begin{document}
 \title{Invariants under deformation of  the actions of topological groups}

 \author{ ANDR\'{E}S   VI\~{N}A}
\address{Departamento de F\'{i}sica. Universidad de Oviedo.   Avda Calvo
 Sotelo.     33007 Oviedo. Spain. }
 \email{vina@uniovi.es}

 \maketitle
\begin{abstract}

Let $\varphi$ and $\varphi'$  be two homotopic actions of the
topological group $G$ on the topological space $X$. To an object
$A$ in the $G$-equivariant derived category $D_{\varphi}(X)$ of
$X$ relative to the action $\varphi$ we associate an object $A'$
of category $D_{\varphi'}(X)$, such that the corresponding
$G$-equivariant compactly supported cohomologies $H_{G,c}(X,\,A)$
and $H_{G,c}(X,\,A')$ are isomorphic.
  When $G$ is a Lie group and $X$ is a subanalytic space, we prove that the
  $G$-equivariant cohomologies $H_{G}(X,\,A)$ and
$H_{G}(X,\,A')$ are also isomorphic.

\end{abstract}

\section {Introduction} \label{S:intro}

\noindent
 Palais and Stewart    proved in \cite{P-S} that the actions of a {\em
compact} Lie group $G$ on a {\em closed smooth manifold} are
rigid; that is, if $\psi$ and $\psi'$ are homotopic $G$-actions on
a closed manifold $N$, then there exists a diffeomorphism $h$ of
$N$, such that $h\circ\psi\circ h^{-1}=\psi'$ (see also \cite[page 191]{G-G-K}). Thus, in
particular, $h$ defines a bijective correspondence between the
fixed point sets of $\psi$ and $\psi'$.


The situation is different, when either the group $G$ is non
compact, or the space on which it acts is not a smooth manifold. A simple example of non rigidity of ${\mathbb R}$-actions is given
in \cite{V3}, where are defined two homotopic ${\mathbb
R}$-actions on $S^1$,  for which there is no a diffeomorphism of
$S^1$ intertwining them.
  In what follows,    we show another example of a non rigid  action of a non compact Lie group
 on ${\mathbb R}P^1$.

 Let $G$ be the Borel subgroup of ${\rm
SL}(2,\,{\mathbb R})$ consisting of the matrices
$$g=\begin{pmatrix} r& y\\
0&r^{-1}
\end{pmatrix},\;\; r\in{\mathbb R}^{\times},\; y\in{\mathbb R}.$$
For $s\in[0,\,1]$, we define the Lie group homomorphism
$\psi^s:G\to{\rm Diff}({\mathbb R}P^1)$  putting
$$\psi^s(g)([a:b])=[ra+syb:r^{-1}b].$$
$\{\psi^s\}_{s\in[0,\,1]}$ is a homotopy between $\psi^0$ and
$\psi^1$.
 The fixed point set of the action $\psi^0$ is $\{[0:1],\,[1:0]\}$,
but the set of fixed points of $\psi^1$ is $\{[1:0]\}$. Thus,
there is no a diffeomorphism of ${\mathbb R}P^1$ which intertwines
$\psi^0$ and $\psi^1$.

 \smallskip

 The preceding examples show that  the  actions of topological groups  do
 admit, in general, non trivial deformations.
The purpose of this note is to   determine  objects
  which are invariants under  deformations of the action of    a topological group $G$ on a topological
 space $X$.
 The underlying general idea we will develop is the following:
 One can expect that homotopic actions on a topological space $X$ give rise to homeomorphic
 homotopy quotients of $X$. On the other hand, the equivariant cohomologies
 are essentially the cohomology of these quotients. Therefore, one can expect
 the existence of an isomorphism between  the equivariant cohomologies of $X$ with coefficients in sheaves which are related through the homotopy.

\smallskip

\noindent
 {\bf Statement of main results.}
 Given a topological space $Y$,
  we write  ${\it Sh}(Y)$ for  the
 abelian category of sheaves of   ${\mathbb C}$-vector spaces on $Y$. By $D^b(Y)$ we denote the
  bounded derived category of ${\it Sh}(Y)$ \cite{Kas-Sch}.

 Associated with the action of a topological group $G$ on the topological space $X$
 there is the homotopy quotient $X_G:=EG\times_G X$,
 and the bounded equivariant
 derived category $D^b_G(X)$ (see \cite{Be-Lu}). Basically, an object $A$ of $D^b_G(X)$
 is a triple $(A_X,\,\bar{ A},\, a)$ consisting of an object ${ A}_X$ of
 the derived category $D^b(X)$, an object $\bar{ A}\in D^b(X_G)$ and $a$ is
  an isomorphism of $D^b(EG\times X)$ between the inverse images of ${A}_X$ and $\bar{A}$
  by the natural maps of this diagram
 $$X\leftarrow EG\times X\rightarrow X_G.$$

The $G$-equivariant cohomology $H_G(X;\,A)$ of $X$ with values in
$A$ is the cohomology of the homotopy quotient
 $X_G$ with coefficients in $\bar A$; i.e., $H(X_G;\,\bar A).$


  By ${\mathcal H}$ we will denote a subgroup of   the  homeomorphism group of $X$.
    Let $\varphi: G\to {\mathcal
   H}$ be a group homomorphism; that is, a $G$-action on $X$.
  From now on,  $D^b_{\varphi}(X)$ stands for
 the equivariant derived category associated with the $G$-action
$\varphi$ on $X$. On the other hand, the homomorphism  $\varphi$
   induces a map $\Phi:BG\to B{\mathcal H}$ between the
   respective classifying spaces.  We write  $X(\varphi)$  for the
   quotient $\Phi^*(E{\mathcal H})\times_{\mathcal H} X$, where
   $\Phi^*(E{\mathcal H})$ is the pullback by $\Phi$ of the
   universal ${\mathcal H}$-bundle $E{\mathcal H}$.

    We will denote by ${\mathfrak
   D}^b(X,\varphi)$ the category whose objects are triples ${\mathcal
   A}=({\mathcal A}_X,\, \bar{\mathcal A},\,\alpha)$, where ${\mathcal
   A}_X\in D^b(X),$ $\bar{\mathcal A}\in D^b(X(\varphi))$ and
   $\alpha$ an isomorphism between the inverse images of ${\mathcal
   A}_X$ and $\bar{\mathcal A}$ by the projections maps in
   $$X\leftarrow \Phi^*(E{\mathcal H})\times X\rightarrow X(\varphi).$$
 Using the universal property of the pullback $\Phi^*(E{\mathcal
H})$, we will construct a functor
   $$\delta:{\mathfrak D}^b(X,\varphi)\to D^b_{\varphi}(X).$$
    Roughly speaking, $\delta$ is a ``forgetful" functor from the ${\mathcal H}$-equivariance to the weaker $G$-equivariance.

   Let $\varphi,\varphi':G\to {\mathcal H}$ be  two  $G$-actions on $X$
   which are homotopic by means of a family of {\em group homomorphisms} from $G$ to ${\mathcal H}$.
   The homotopy gives rise to
    an equivalence $\xi^*$ between the categories
    ${\mathfrak D}^b(X,\varphi')$ and ${\mathfrak
    D}^b(X,\varphi)$.

The $G$-equivariant (relative to the action $\varphi$) cohomology
of $X$ with compact supports and coefficients in $A\in
D^b_{\varphi}(X)$ will be denoted by $H_{\varphi,c}(X;\,A)$. We
will prove the following theorem.

 \begin{Thm}\label{Thm:EquivCohom}
Let $G$ be a topological group   and  $\varphi,\varphi':G\to
{\mathcal H}$ be continuous group homomorphisms that are homotopic
through a family
 of group homomorphisms from $G$ to ${\mathcal H}$.
 If $X$ and $G$  are locally compact,
 then the cohomologies  $H_{\varphi,c}(X;\,A)$ and $H_{\varphi',c}(X;\,A')$ are isomorphic
$$H_{\varphi,c}(X;\,A)\simeq H_{\varphi',c}(X;\,A'),$$
  where
$A=\delta({\mathcal A})$, $A'=\delta'(\xi^{*})^{-1}({\mathcal A})$
and ${\mathcal A}$ is any object of  ${\mathfrak D}^b(X,
\varphi)$.
 \end{Thm}



In the case when $G$ is a   Lie group and $X$ is a subanalytic
space, we prove
the existence of an isomorphism between the $G$-equivariant
cohomologies $H_{\varphi}(X;\,A)$ and $H_{\varphi'}(X;\,A')$.
 More precisely,
 \begin{Thm}\label{Thm:EquivCohom2}
 Let $\varphi,\varphi':G\to {\mathcal H}$ be group actions on $X$, such that
\begin{enumerate}
\item $\varphi$ and $\varphi'$  are homotopic through a family
 of group homomorphisms from $G$ to ${\mathcal H}$,
 \item $X$ is  a subanalytic space,
\item $G$ is a Lie group with a finite number
of connected components.
\end{enumerate}
 Then the cohomologies   $ H_{\varphi}(X;\,A)$ and  $ H_{\varphi'}(X;\,A')$  are isomorphic
 $$ H_{\varphi'}(X;\,A')\simeq H_{\varphi}(X;\,A),$$
where $A=\delta({\mathcal A})$,
$A'=\delta'(\xi^{*})^{-1}({\mathcal A})$ and ${\mathcal A}$ is any
object of ${\mathfrak D}^b(X,\varphi)$.
 \end{Thm}


  In the case that the group $G$ is a compact torus $T$ and $X$ is a compactifiable
  space (see \cite[Sections 3 and 6]{G-K-M}),
  we will prove
  the following theorem.

\begin{Thm}\label{Cohom=}
 Let $X$ be a compactifiable space  and $\varphi,\,\varphi':T\to {\mathcal H}$ be continuous group homomorphisms
 homotopic through a family of group homomorphisms.
   Denoting by    $F$ and $F'$   the  fixed point sets of the corresponding actions on $X$, then
  $$ H(F;\,{\mathbb C})\simeq H(F';\,{\mathbb C}).$$
 \end{Thm}

\smallskip

\section{Equivariant derived categories}

\subsection{Homotopy quotients}\label{SubHomotopy quotients}

 In this subsection, we define the spaces
$X(\varphi)(n)$, which are ``approximations" to the space
$X(\varphi)$ introduced in Section \ref{S:intro}.

Let $X$ be a Hausdorff locally compact topological space and $G$ a
topological group. As we said, ${\mathcal H}$ stands for
a topological subgroup of the homeomorphism group of $X$, endowed
with the compact open topology.  By $\varphi:G\to {\mathcal H}$ we
denote a continuous group homomorphism.

Given a positive integer $n$, let $EG(n)\to BG(n)$ be the
 approximation to the universal $G$-bundle  given by   the $n$-fold joint
$\overbrace{G\ast\cdots\ast G}^n$ in the Milnor's construction
(see \cite{Mi}, \cite[page 53]{Hus}). Following Husemoller, we
write $\langle g,t\rangle$ for the elements of $EG(n)$, where
$\langle g,t\rangle$ is the
 sequence $\langle t_1g_1,\dots ,t_ng_n\rangle$ satisfying the properties detailed  in
 \cite{Hus}.
As it is well-known, the topological space $EG(n)$ is
$\kappa(n)$-acyclic, where $\kappa(n)$ is a number that tends to
infinity when $n$ grows.

It is easy to check in the Milnor's construction that the group
homomorphism $\varphi$ induces a continuous map $\phi_n:EG(n)\to
E{\mathcal H}(n)$, satisfying $\phi_n(eg)=\phi_n(e)\varphi(g)$,
for all $e\in EG(n)$ and all $g\in G$.
In fact, in the above notation the map $\phi_n$ is given by
 \begin{equation}\label{phnlangle}
 \phi_n(\langle
g,t\rangle)=\langle\varphi(g),t\rangle.
 \end{equation}

 In turn,
$\phi_n$ induces a continuous map $\Phi_n:BG(n)\to B{\mathcal
H}(n)$, between the corresponding approximations to the
classifying spaces. The map $\phi_n$ factors through
$\Phi^*_n(E{\mathcal H}(n))$, the pullback of $E{\mathcal H}(n)$
by the map $\Phi_n$,
\begin{equation}\label{tau(n)}
 \xymatrix{
EG(n) \ar[dr]_-{\phi_n}\ar[rr]^-{\tau_n} &&\Phi_n^*(E{\mathcal H}(n))\ar[dl]\\
&E{\mathcal H}(n)      }
\end{equation}

The $G$-action on $X$ determines an associated bundle to $EG(n)$
with fiber $X$, that will be denoted
$X_{\varphi}(n):=EG(n)\times_{\varphi}X$. We have also the
corresponding homotopy quotient $X_{\mathcal H}(n):=E{\mathcal H}(n)\times_{\mathcal H}X$.
 The map between the associated bundles  $X_{\varphi}(n)$ and $X_{{\mathcal H}}(n)$ induced
by $\phi_n$,   will be denoted by $f_n$.
 The map $f_n$ factors through $X(\varphi)(n):=\Phi_n^*(E{\mathcal H}(n))\times_{\mathcal H} X$;
$f_n=\nu_n\circ\rho_n$.
$$\xymatrix{
X_{\varphi}(n) \ar[dr]_-{f_n}\ar[rr]^-{\rho_n} &&  X(\varphi)(n)\ar[dl]^-{\nu_n}\\
&X_{\mathcal H}(n)      }
$$

 In  this  subsection, we assume that the value
of $n$ has been fixed  and we will {\em omit} the $n$'s in
$BG(n)$, $EG(n)$, $X_{\varphi}(n)$, $X(\varphi)(n)$, $\phi_n$,
$\rho_n$, $f_n$, etc.

Some objects and arrows we have introduced can be organized in the
 following commutative diagram in the category of topological
spaces
\begin{equation*}\tag{cd(n)}
 \xymatrix{ X  \ar@<0.5ex>[d]^{1} &    EG\times X\ar@<0.0ex>[l]_{p} \ar@<0.0ex>[r]^{q}
\ar@<0.5ex>[d]^{\tau\times 1}
 & X_{\varphi} \ar[d]^{\rho}  \\
X & \Phi^*(E{{\mathcal H}}) \times
X\ar@<0.0ex>[r]^-{Q}\ar@<0.0ex>[l]_-{P} & X(\varphi)
  \; , }
\end{equation*}
where $p$ and $P$ are the corresponding projection maps and $q$
and $Q$ the quotient maps.

 In the following   commutative diagram  are
involved the bundles with fiber  $X$ that we have mentioned

\begin{equation} \label{doublesquare}
  \xymatrix{ X_{\varphi}\ar@(ur,ul)[rr]^{f}\ar@<0.0ex>[r]^-{\rho}  \ar@<0.5ex>[d]^{\pi} &
 X(\varphi) \ar@<0.0ex>[r]^-{\nu}
\ar@<0.5ex>[d]^{\bar\pi}
 & X_{\mathcal H} \ar[d]^{\pi^{\mathcal H}}  \\
BG\ar@<0.0ex>[r]^{1} & BG\ar@<0.0ex>[r]^{\Phi} & B{\mathcal H}
  \; . }
\end{equation}

\begin{Prop}\label{compactkernel} The map $\rho:X_{\varphi} \to X(\varphi)$ is an isomorphism of
bundles over the identity of $BG$.
\end{Prop}

{\it Proof.} The mapping $\rho$ is given by $\rho(q(e,x))=Q(\phi(e),x)$, for $e\in BG$ and $x\in X$.
 Fixed an
 element $e^0$ in the fibre of $b$ in the fibration $EG\to BG$,
 then $\pi^{-1}(b)=\{ q(e^0,\,x) \,|\, x\in X\}$ and $\pi^{-1}(b)$ can be parametrized by the points of $X$. An analogous
 parametrization of $\bar\pi^{-1}(b)$ is defined by choosing an
 element $e'$ in the fibre over $b$ in $\Phi^*(E{\mathcal H})\to BG$;
 in particular, we can choose $e'=\phi(e^0)$. With these choices the
 map $\pi^{-1}(b)\to\bar\pi^{-1}(b)$, restriction of $\rho$, is
 simply the map $x\mapsto x$. That is,
$\rho$ induces homeomorphisms between $\pi^{-1}(b)$ and $\bar\pi^{-1}(b)$,
 for all $b\in BG$.

 For $i\in \{1,\dots,n\}$, let $t_i$ be the natural map $\langle g,t\rangle\in EG\mapsto t_i\in [0,\,1]$ (see  \cite{Hus}).
 The family
 \begin{equation}\label{CoveringBG}
 \{V_i:= \Hat\pi(t_i^{-1}(0,\,1 ])\}_{i=1,\dots,n},
  \end{equation}
  $\Hat\pi$ being the projection $EG\to BG$,    is an open covering of $BG$.

  On $V_i$ we consider the local section $\chi_i$ of $EG$ defined by
  $$\chi_i([g,\,t])=\langle g,\,t\rangle g_i^{-1}=\langle t_1g_1g_i^{-1},\dots, t_ng_ng_i^{-1}\rangle,$$
  where $[g,\,t]=\Hat\pi(\langle g,\,t\rangle)$.
   The equality $[g,t]=[g',t']$ implies the existence of an element
 $b\in G$ such that $g'_j=g_jb$, for all $j$ such that $t_j\ne 0$.
 If  $[g,t]=[g',t'] \in V_i$, then
   $b=g_i^{-1}g'_i$; thus, $ g_j g_i^{-1}=g'_jg_i'^{-1}$, and
    the section $\chi_i$ is well defined.
 On $V_i\cap V_j$, $\chi_j=\chi_i\cdot m_{ij}$, with
$m_{ij}([g,\,t])=g_ig_j^{-1}.$

    Analogously, $\{\tilde V_i\}_i$, with
    $$\tilde V_i=\{[h,t]\,|\, t_i\ne 0  \}\subset B{\mathcal H},$$
    is an open covering of $B{\mathcal H}$. On each $\tilde V_i$ one can define the corresponding section $\tilde\chi_i$
    of $E{\mathcal H}.$ The respective transition functions are $\tilde m_{ij}([h,\,t])=h_ih_j^{-1}$.

    Since the map $\Phi:BG\to B{\mathcal H}$ is defined by $\Phi([g,\,t])=[\varphi(g),\,t]$,
    the inverse image $\Phi^{-1}(\tilde V_i)=V_i$, and the transition
     functions of $\Phi^*(E{\mathcal H})$ associated with the above trivializations are the maps
     $$[g,\,t]\in V_i\cap V_j\mapsto \varphi(g_i)\varphi (g_j^{-1})\in {\mathcal H}.$$

    On other hand, $X_{\varphi}$ is the bundle associated to $EG$ through the representation $\varphi$ of $G$ in the homeomorphisms
   of $X$. Therefore, its transition functions are the functions $[g,\,t]\to\varphi(g_ig_j^{-1})$. Thus, the bundles  $X_{\varphi}$ and $X(\varphi)$
   are isomorphic, since they have  homeomorphic fibers and
    the same transition functions with respect to the families of local sections $\{\chi_i\}_i$ and $\{\tilde\chi_i\}_i$.

   \qed

\begin{Cor}\label{CorIndependenceH}
If $\tilde{\mathcal H}$ is a subgroup of ${\mathcal H}$ such that
$\varphi$ can be expressed as the composition
$i\circ\tilde\varphi$, with $i:\tilde{\mathcal
H}\hookrightarrow{\mathcal H}$, then $ X(\varphi)$ and $
X(\tilde\varphi)$ are isomorphic bundles.
\end{Cor}

  If $X$ is a  locally compact  topological space and
 $G$ is a  locally
 compact topological group, then
   the spaces $ X(\varphi)\simeq X_{\varphi} $ and $BG$
 are locally compact, as well, and we can consider functors direct
  image with compact support between the categories of sheaves
 on these spaces.

By Proposition \ref{compactkernel},  the left square in (\ref{doublesquare}) is obviously cartesian. Then we have the
 ``base change" formulas given in the following proposition.

\begin{Prop}\label{basechange}
 Let assume that
\begin{enumerate}
\item  $X$ is a locally compact topological space,
 \item $G$ is   locally
 compact topological group.
 \end{enumerate}
 Then there is a
canonical isomorphism between the functors $R\pi_!\rho^*$ and
$R\bar\pi_!$
 \begin{equation}\label{basechange1}
 R\pi_!\rho^*\simeq R\bar\pi_!:D^b(X(\varphi))\to D^b(BG).
\end{equation}
And a canonical isomorphism between the functors
  \begin{equation}\label{basechange2}
  R\pi_*\rho^*\simeq R\bar\pi_*: D^b(X(\varphi))\to D^b(BG).
 \end{equation}
 \end{Prop}

\subsection{The space $X(\varphi)$} \label{ThespacebarX}
 Next, we define the space $X(\varphi)$ as the direct limit of the
$X(\varphi)(n)$'s. We insert the label $n$ in the notations. Let
us consider the filtrations associated with the group $G$
 $$\dots\subset EG(n)\subset EG(n+1)\subset\dots \;\; \hbox{\rm and}\;\; \dots\subset BG(n)\subset
 BG(n+1)\subset\dots$$
 and the ones corresponding to the group ${\mathcal H}$. If $\varphi:G\to{\mathcal H}$ is a $G$-action on $X$,
 the maps
 \begin{equation}\label{phin}
  \Phi_n:BG(n)\to B{\mathcal H}(n)    \;\; \hbox{and}\;\;  \tau_n:EG(n)\to \Phi_n^*(E{\mathcal H}(n))
 \end{equation}
 are
 compatible with the filtrations. That is, for $n<m$ we denote
 $$j_{mn}:EG(n)\hookrightarrow  EG(m),\;\; \hbox{and}\;\; i_{mn}:\Phi_n^*(E{\mathcal H}(n))\hookrightarrow
 \Phi_m^*(E{\mathcal H}(m)),$$
  then $\tau_mj_{mn}=i_{mn}\tau_n$.
Similar relations are valid for the maps $\Phi_n$ and the
 inclusions $BG(n)\hookrightarrow BG(m)$.
On the other hand,    $j_{mn}$ induces a map $\psi_{mn}:X_{\varphi}(n)\to X_{\varphi}(m)$,
  and in turn $i_{mn}$ induces   $\bar\psi_{mn}:X(\varphi)(n)\to
 X(\varphi)(m)$. In summary, the following diagrams are commmutative
 $$\begin{matrix}
      \xymatrix{   \Phi_n^*(E{\mathcal H}(n))      \ar@<0.0ex>[r]^{i_{mn}}
 &  \Phi_n^*(E{\mathcal H}(m)) \\
EG(n)\ar[u]^{\tau_n} \ar@<0.0ex>[r]^{j_{mn}}   & EG(m)\ar[u]_{\tau_m}   }   \;\;\;    \;\;\;      \;\;\;                 &

  \xymatrix{ X(\varphi)(n)    \ar@<0.0ex>[r]^{\bar\psi_{mn}}
 & X(\varphi)(m)   \\
 X_{\varphi}(n)\ar[u]^{\rho_n}  \ar@<0.0ex>[r]^{\psi_{mn}}   & X_{\varphi}(m)\ar[u]_{\rho_m}  }
 \end{matrix}
 $$

For  later references,  we copy left square of
(\ref{doublesquare})   showing the dependence of $n$ of its
elements.

\begin{equation*}\tag{sq(n)}
\xymatrix{ X_{\varphi}(n) \ar@<0.0ex>[r]^{\rho_n}
\ar@<0.5ex>[d]^{\pi_n}
 & X(\varphi)(n) \ar[d]^{\bar\pi_n}  \\
BG(n) \ar@<0.0ex>[r]^{1}   & BG(n) \; . }
\end{equation*}

 The maps $\psi_{mn}$, $\bar\psi_{mn}$ and  the inclusion
$BG(n)\subset BG(m)$  connect the vertices of diagram (sq(n))
 with the corresponding vertices of (sq(m)), and the resulting cube is
 a commutative diagram. Taking the inductive limit, we obtain
 a homeomorphism
 $$\rho:X_{\varphi}:=\underset{\rightarrow}{\lim}\,X_{\varphi}(n)\longrightarrow
 X(\varphi):=\underset{\rightarrow}{\lim}\,X(\varphi)(n),$$
  which is an isomorphism over the identity between the bundles $X_{\varphi}\overset{\pi}{\to} BG$ and
  $X({\varphi}) \overset{\bar\pi}{\to} BG.$



\smallskip

\subsection{The category ${\mathfrak D}^b(X,\varphi)$.} \label{The category}

In this subsection, we assume over the space  $X$ and the group $G$ the hypotheses of Proposition \ref{basechange}.
Then, as we said, the
 spaces $X_{\varphi}(n)$ and  $X(\varphi)(n)$
 are also locally compact.

 To avoid the use of non locally compact topological spaces in the definition of  ${\mathfrak D}^b(X,\varphi)$, we will construct this category as a ``direct limit" of a sequence of categories, following a procedure similar to the one employed for  the definition of the $G$-equivariant category $D^b_G(X)$  in \cite[page 27]{Be-Lu}.

\smallskip
\noindent
 {\it The category ${\mathfrak D}^b(X,\varphi)$.} Let
${\mathcal A}_n:=({\mathcal A}_X,\bar{\mathcal A}_n,\alpha_n)$ be
a triple consisting of
$${\mathcal A}_X\in D^b(X),\;\;\;\;\bar{\mathcal
A}_n\in D^b(X(\varphi)(n))$$
 and $\alpha_n$ an isomorphism in
$D^b\big(\Phi^*_n(E{\mathcal H}(n))\times X\big)$ from
$P_n^*({\mathcal A}_X)$ to $Q^*_n(\bar{\mathcal A}_n)$ (see
diagram (cd(n))). That is, with
  the notation introduced in  \cite[page 17]{Be-Lu},
the triple ${\mathcal A}_n$  defines
 an object  of the   category $D^b_{ \mathcal H }(X,\, \Phi^*_n(E{\mathcal H}(n)) \times X)$.

Let ${\mathcal A}:=\{{\mathcal A}_n\}_{n=1,2,\dots}$ be a sequence
of triples, such that each of which satisfies the properties
stated above. Furthermore, we assume that the following
compatibility conditions hold
 \begin{equation}\label{compti}
 \bar\psi^*_{mn}(\bar{\mathcal A}_m)= \bar{\mathcal A}_n,\;\;\;(i_{mn}\times 1)^*(\alpha_m)=\alpha_n,
  \end{equation}
 for $n<m$.

The above sequence ${\mathcal A}$ does not defines an object of
the category $D^b_{\mathcal H}(X)$, since the limit of the
$\Phi_n^*(E{\mathcal H}(n))$'s is not contractible. However,
we can construct a new category, ${\mathfrak D}^b(X,\varphi)$,
whose objects are   the ${\mathcal A}$'s.

 Given two objects ${\mathcal A}=\{({\mathcal A}_X, \bar{\mathcal A}_n,\alpha_n)\}_n$ and
 ${\mathcal B}=\{({\mathcal B}_X, \bar{\mathcal B}_n,\beta_n)\}_n$ in the
 category ${\mathfrak D}^b(X,\varphi)$,
  a morphism   $\sigma:{\mathcal A}\to{\mathcal B}$ in ${\mathfrak D}^b(X,\varphi)$ is a sequence of pairs
  $\{(\sigma_X,\bar\sigma_n)\}_n$, with
 $\sigma_X:{\mathcal A}_X\to{\mathcal B}_X$   and $\bar\sigma_n:\bar{\mathcal A}_n\to\bar{\mathcal B}_n$, such that
 $$\beta_n \circ P_n^*(\sigma_X)=Q_n^*(\bar\sigma_n)\circ\alpha_n,$$
  and for $n<m$
 they satisfy the
 compatibility condition
 \begin{equation}\label{compatibility}
 \bar\psi^*_{mn}(\bar\sigma_m)=\bar\sigma_n.
 \end{equation}

\smallskip
\noindent
 {\it The functor $\delta$.}  On the other hand,
 ${\mathcal A}_n$
gives rise to the triple $A_n:=(A_X,\bar A_n,a_n)$, where
\begin{equation}\label{f*}
A_X={\mathcal A}_X\in D^b(X),\;\;\bar A_n=\rho_n^*(\bar{\mathcal
A}_n)\in D^b(X_{\varphi}(n)),
\end{equation}
and
\begin{equation}\label{a_n}
 a_n:=(\tau_n\times 1)^*(\alpha_n)
 \end{equation}
  is an isomorphism in $D^b(EG(n)\times X)$ between
$p_n^*(A_X)\to q_n^*(\bar A)$ (see cd(n)). In other words, $A_n$
is an object of $D^b_{\varphi}(X,\; EG(n)\times X)$,
  when  we
consider in $X$ the $G$-action defined by $\varphi$.

Given the object ${\mathcal A}=\{\mathcal A_n\}$,  each ${\mathcal
A}_n$ gives rise to a triple $A_n=(A_X,\bar A_n,a_n)$. Since
diagrams cd(n) and cd(m) together with the maps between them
 determined by $j_{mn},$ $i_{mn}$, $\psi_{mn}$ and $\tilde
 \psi_{mn}$ form a commutative cuboid, then
from (\ref{compti}), it follows
 \begin{equation}\label{compti1}
  \psi^*_{mn}(\bar{ A}_m)= \bar{ A}_n,\;\;\;(j_{mn}\times 1)^*(a_m)=a_n.
  \end{equation}

As the $G$-space $EG(n)\times X$ defines a resolution
$\kappa(n)$-acyclic of $X$,  the sequence $A=\{A_n\}_n$ can be
considered as an object  of the equivariant derived category
$D^b_{\varphi}(X)$ \cite[page 27]{Be-Lu}.

On the other hand, given   $\sigma:{\mathcal A}\to{\mathcal B}$,
we define $s_X:=\sigma_X$ and $\bar s_n:=\rho_n^*(\bar\sigma_n)$.
 From
 (\ref{a_n}) together with the commutativity of (cd(n)), it follows $q_n^*(\bar s_n)\circ a_n=b_n\circ p_n^*(s_X)$.
  By (\ref{compatibility}),
 $$\psi^*_{mn}(\bar s_m)=\bar s_n.$$
 Thus, the family  $s=\{(s_X,\bar s_n)\}_n$
 is a morphism in $D^b_{\varphi}(X)$ from $A$ to $B$.
So, we have the functor
 \begin{equation}\label{Funtornu}
 \delta:{\mathfrak D}^b(X,\varphi)\to D^b_{\varphi}(X);\;\;{\mathcal A}\mapsto A,\;\;
 \sigma\mapsto s.
  \end{equation}





\smallskip

{\it Examples}.  With the
constant sheaves ${\mathbb C}_X$ and
  ${\mathbb C}_{X(\varphi)(n)}$, we can construct the triple
  $${\mathcal C}_n:=({\mathbb C}_X,\,{\mathbb C}_{X(\varphi)(n)},\, {\rm Id})\in D^b_{
   \mathcal H }\big(X,\,\Phi_n^*(E{\mathcal H}(n))\times X\big),$$
since the constant sheaf on a space $Z$ is defined through the inverse
image by $\mu:Z\to {\rm pt}$,  the constant  map to a point. Then,
the object $\delta({\mathcal C})\in D^b_{\varphi}(X)$
 is defined by the
sequence $\{C_n\}$, with
$$C_n:=({\mathbb C}_X,\,{\mathbb C}_{X_{\varphi}(n)},\,
 {\rm Id})\in D^b_{\varphi}(X,\; EG(n)\times X) ,$$
 again by the above mentioned
  property of the constant sheaf.

\smallskip

  Let us assume that $X$ is a topological space with finite cohomological dimension. The dualizing object $D_X$ on $X$ is defined as the inverse image with compact support  of ${\mathbb C}_{\rm pt}$ by the constant map to a point. Since $Q^*\mu^!=\mu^!$ and
  $P^*\mu^!=\mu^!$, we can construct the following object of
   the category $D^b_{ \mathcal H }\big(X,\,\Phi_n^*(E{\mathcal H}(n))\times X\big)$
  \begin{equation}\label{Ddualizing}
  {\mathcal D}_n=(D_X,\,D_{X(\varphi)(n)},\, {\rm Id})\in D^b_{
   \mathcal H }\big(X,\,\Phi_n^*(E{\mathcal H}(n))\times X\big),
   \end{equation}
   which gives rise to the object $\delta({\mathcal D})\in D^b_{\varphi}(X).$

\smallskip

{\it Remark.} If $\tilde {\mathcal H}$ is a subgroup of ${\mathcal H}$ containing the image of $\varphi$, as in Corollary \ref{CorIndependenceH}, the inclusion
$E\tilde{\mathcal H}\subset E{\mathcal H}$ induces
a factorization of $\rho_n:X_{\varphi}(n)\to X(\varphi)(n)$
through $ X(\tilde\varphi)(n)$,  which in turn gives rise to the
following factorization of the functor $\delta$
\begin{equation}\label{factorizationrho}
 \xymatrix{
{\mathfrak D}^b(X,\varphi)\ \ar[dr]_{\delta}  \ar[rr] && {\mathfrak D}^b( X,\tilde\varphi) \ar[dl]^{\tilde\delta}   \\
&    D^b_{\varphi}(X)  \;.  }
\end{equation}

 According to Corollary \ref{CorIndependenceH}, the spaces $X(\varphi)$ and   $X(\tilde\varphi)$ are homeomorphic, but
 the categories ${\mathfrak D}^b(X,\varphi)$ and ${\mathfrak D}^b( X,\tilde\varphi)$ are not necessarily equivalent.


\section{Deformation of group actions}\label{SectDeformation of group actions}
 Before to consider deformations of $G$-group actions on the space
 $X$, we will treat   with a slighter general
 situation.

Given  $\varphi,\,\varphi':G\to {\mathcal H}$  two
  $G$-actions {\em not} necessarily homotopic. Let us assume that
  for each $n$ there are bundle isomorphisms over the identity of
  $BG(n)$
  $$\omega_n: \Phi^*_n(E{\mathcal H}(n))\to \Phi'^*_n(E{\mathcal
H}(n))\;\;\hbox{and}\;\; \vartheta_n: X(\varphi)(n) \to
X(\varphi')(n),$$
   satisfying  the compatibility conditions
   \begin{equation}\label{itemizef}
   \vartheta_n\circ Q_n=Q'_n\circ(\omega_n\times 1), \;\;\;
  i'_{mn}\circ\omega_n=\omega_m\circ i_{mn},\;\;\;\bar\psi'_{mn}\circ\vartheta_n=\vartheta_m\circ\bar\psi_{mn},
 \;\;\hbox{for}\; n<m.
 \end{equation}
  For the sake of brevity,  we say that the family
  $\vartheta=\{\vartheta_n\}$ is a $(\varphi,\,\varphi')$-cohomological
 morphism. (The reason for the adjective cohomological will be
 apparent below).

That is, for each $n$, we have the following commutative diagram
of fiber bundles on $BG(n)$, that in turn is consistent with the
one corresponding to the number $m$.
\begin{equation}\label{DiagramIV}
\xymatrix{
  {} &  EG(n)\times X\ar@{-->}[rr]^{q'_n}\ar@{-->}[dd]^(0.3){\tau'_n\times 1}  & {}  &  X_{\varphi'}(n) \ar@{<-->}[dd]^(.6){\rho'_n} \\
  EG(n)\times X \ar[rr]^(.3){q_n}\ar[dd]^(0.3){\tau_n\times 1}  & {}  & X_{\varphi}(n) \ar@{<->}[dd]_(.3){\rho_n} & {}  \\
  {} & {\Phi'_n}^*(E{\mathcal H}(n))\times X \ar@{-->}[rr]^(.7){Q'_n} &  {} &   X({\varphi'_n})  \\
  \Phi_n^*(E{\mathcal H}(n))\times X \ar[rr]^{Q_n}\ar@{<..>}[ru]_{\omega_n\times 1} & {} &  X({\varphi}_n) \ar@{<..>}[ur]_(0.5){\vartheta_n}  & {} \;.
         }
 \end{equation}

\begin{Prop}\label{isomorFunctors}
  If  the space $X$ and the group $G$   are locally
  compact and $\vartheta$ is a $(\varphi,\,\varphi')$-cohomological morphism, then
  there exist the following canonical isomorphisms of functors
  $$R(\bar\pi'_n)_!\simeq
  R(\bar\pi_n)_!\vartheta_n^*: D^b(X(\varphi')(n))\to D^b(BG(n))$$
 $$R(\bar\pi'_n)_*\simeq
  R(\bar\pi_n)_*\vartheta_n^* : D^b(X(\varphi')(n))\to D^b(BG(n)).$$
  \end{Prop}
{\it Proof.} Since $\vartheta_n$ is a homeomorphism, the diagram
 \begin{equation}\label{DiagramII}
\xymatrix{ X(\varphi)(n) \ar@<0.0ex>[r]^{\vartheta_n}
\ar@<0.5ex>[d]^{\bar\pi_n}
 &   X(\varphi')(n) \ar[d]^{\bar\pi'_n}  \\
BG(n) \ar@<0.0ex>[r]^{1}   & BG(n) \; . }
\end{equation}
  is a
  cartesian square, obviously. Formula (2.4b) in \cite{G-Mc} applied to this square gives the first isomorphism.
  On the other hand, the square consisting of the maps $\vartheta_n^{-1}$, $\bar\pi_n$, $1$ and $\bar\pi'_n$
  is also cartesian. Applying (2.5a) of \cite{G-Mc} to those cartesian squares, one obtains the second isomorphism.

 \qed

Taking the inductive limits as $n$ goes to $\infty$ in
(\ref{DiagramII}), one  prove that the diagram obtained from
(\ref{DiagramII}) deleting the label $n$ is also commutative.

 \begin{Prop}\label{PropItemize}
  The families $\omega$  and $\vartheta$ define an isomorphism $\vartheta^*$ from the  category
${\mathfrak D}^b(X,\varphi')$ to ${\mathfrak D}^b(X,\varphi)$.
\end{Prop}

{\it Proof.}
  The inverse image functors induced by the family of homeomorphisms  $\omega_n$ and  $\vartheta_n$
   determine the isomorphism.
   More precisely, the functor
   $(\vartheta^{-1})^*:{\mathfrak D}^b(X,\varphi)\to {\mathfrak D}^b(X,\varphi')$ defined by
   $$(\vartheta^{-1})^*({\mathcal A})=\big\{\big({\mathcal A}_X,\,(\vartheta_n^{-1})^*(\bar{\mathcal A}_n),
   \,(\omega_n^{-1}\times 1)^*(\alpha_n)\big)   \big\}_n$$
     gives the isomorphism between the categories.
 \qed

\begin{Prop}\label{PropEnumerate} Under the hypotheses of
Proposition \ref{isomorFunctors}, given
 ${\mathcal A}\in{\mathfrak D}^b(X,\varphi)$,
 let
 $$\{(A_X,\bar A_n, a_n)\}=\delta({\mathcal A})\in D^b_{\varphi}(X) \;\;\;
 \hbox{and} \;\;\; \{(A_X,\bar A'_n, a'_n)\}=\delta'(\vartheta^{-1})^*({\mathcal A})\in D^b_{\varphi'}(X),$$ then
\begin{itemize}
 \item The objects $R(\pi_{n})_!(\bar A_n)$ and
$ R(\pi'_n)_{!}(\bar A'_n)$  of $D^b(BG(n))$  are canonically
isomorphic.
 \item   $R(\pi_{n})_*(\bar A_n)\simeq  R(\pi'_n)_{*}(\bar A'_n)$
canonically isomorphic.
 \end{itemize}
\end{Prop}
 {\it Proof.} From Proposition \ref{basechange} together with (\ref{f*}),
  it follows
  $$ R(\pi_{n})_!(\bar A_n)=R(\pi_{n})_!\rho_n^*(\bar{\mathcal A}_n)
\simeq R(\bar\pi_n)_{!}(\bar{\mathcal
  A}_n),$$
  where   the isomorphism is canonical.

   Analogously, $R(\pi'_n)_{!}(\bar A'_n)\simeq R(\bar\pi'_n)_{!}(\vartheta_n^{-1})^*(\bar{\mathcal A}_n).$ Using
 Proposition \ref{isomorFunctors}, we conclude
 $$ R(\pi_{n})_!(\bar A_n)\simeq R(\pi'_n)_{!}(\bar A'_n).$$

 Similarly, from Proposition \ref{isomorFunctors} and Proposition
\ref{basechange}, it follows our second claim.

\qed

As above, we  denote by $\mu$ the  constant map $\mu:X\to {\rm
pt}$. Then  we have the following commutative diagram
 \begin{equation}\label{mu}
 \xymatrix{ X  \ar@<0.5ex>[d]^{\mu} &    EG(n)\times X\ar@<0.0ex>[l]_-{p_n} \ar@<0.0ex>[r]^{q_n}
\ar@<0.5ex>[d]^{\epsilon}
 & X_{\varphi}(n) \ar[d]^{\mu_G}  \\
{\rm pt}& EG(n)\ar@<0.0ex>[r]^-{\Hat\pi} \ar@<0.0ex>[l]_{u}&
BG(n)={\rm pt}_{G}(n)
  \; ,}
 \end{equation}
where $\mu_G$ is the map induced between the homotopy quotients;
that is, $\mu_G=\pi_n$ (see (sq(n)).

Let $A=\{(A_X,\bar A_n,a_n)\}$ be an object of $D_{\varphi}^b(X)$,
then
 formula (2.4b) of \cite{G-Mc} applied to the squares of (\ref{mu})
gives the isomorphisms
\begin{equation}\label{changeGoreski}
u^*R\mu_!(A_X)\simeq R\epsilon_!p_n^*(A_X),\;\hbox{and}\,\;
R\epsilon_!q_n^*(\bar A_n)\simeq \Hat\pi^*R (\pi_{n})_!(\bar A_n).
 \end{equation}
 The
object $R\mu_!(A)\in D_G^b({\rm pt})$ is defined by the family of
triples
$$\{\big(R\mu_!(A_X),\,R (\pi_{n})_!(\bar A_n),\,c_n\big)\},$$
 where $c_n$ is given by
the composition of the canonical isomorphisms
(\ref{changeGoreski}) with
$$R\epsilon_!(a_n):R\epsilon_!p_n^*(A_X)\to R\epsilon_!q_n^*(\bar
A_n).$$
 Thus, we have defined a functor $D^b_{\varphi}(X)\to D^b_G({\rm pt})$, which will be denoted by $R\mu_!$.

 From Proposition \ref{PropEnumerate}, it follows the following proposition.
\begin{Prop}\label{Prop;IsomrphismLambda}
 If the group $G$ and the space $X$ are locally compact and $\vartheta$ is
 a $(\varphi,\,\varphi')$-cohomological morphism,
 then the  isomorphism $\vartheta^*:{\mathfrak D}^b(
X,\varphi')\to {\mathfrak D}^b(X,\varphi)$ makes commutative, up
to canonical isomorphism, the following diagram
\begin{equation}\label{Diagram}
 \xymatrix{
 {} & D^b_G({\rm pt})   &{}\\
  D^b_{\varphi'}(X)\ar@<0.0ex>[ur]^{R\mu'_!}
      & {} & D^b_{\varphi}(X)  \ar@<0.0ex>[ul]_{R\mu_!}         \\
  {\mathfrak D}^b( X,\varphi')\ar@<0.0ex>[rr]^{\vartheta^*}\ar@<0.0ex>[u]^{\delta'} & {} &
   {\mathfrak D}^b(X,\varphi)\ar@<0.0ex>[u]_{\delta}\; ,           }
    \end{equation}

where ${R\mu_!}$ and ${R\mu'_!}$ are the functors induced by the
constant map $\mu:X\to{\rm pt}$.
  \end{Prop}

As we said,  we denote by  $H_{\varphi, c}(X;\,A)$ the equivariant
(with respect to the action $\varphi$) cohomology of $X$ with
compact supports
 and coefficients in $A\in D^b_{\varphi}(X)$. By definition   $H_{\varphi,
c}(X; \,A)= H(R\mu_!A)$
\cite[page 115]{Be-Lu}. A direct consequence of above Proposition
is the following corollary.
 \begin{Cor}\label{CorIguality}
Under the hypotheses  of Proposition \ref{Prop;IsomrphismLambda},
 if
${\mathcal A}\in{\mathfrak D}^b(X,\,\varphi)$, then
$$H_{\varphi,c}(X; \,\delta(\mathcal A))\simeq H_{\varphi,c}(X;
\,\delta'(\theta^{-1})^*(\mathcal A)).$$
 \end{Cor}

To introduce the functor
direct image  in the context of equivariant categories are
necessary additional hypotheses on the topological space $X$ and
on the group $G$, because the formula of change basis for the
functor direct image holds only when the change is smooth
\cite[page 10]{Be-Lu}.

   We will assume that $X$ is  a subanalytic space.
 With respect to   $G$, we suppose that it is  a group  such
that the spaces $EG(n)$  are homeomorphic to finite dimensional
manifolds. This property is guaranteed if $G$ is a Lie group with
a finite number of connected components
 (see \cite[pages 34-35]{Be-Lu}).

 As in the case of the direct image with
 compact supports, we consider diagram (\ref{mu}) but assuming that $EG(n)$ is a smooth manifold.
  Given $A=\{(A_X,\,\bar A_n,\,a_n)\}\in D_G^b(X)$, since the space $EG(n)$ is
 smooth  we have the following isomorphisms (see \cite[page 13]{Be-Lu},\cite[page 38]{G-K-M})
$$u^*R\mu_*(A_X)\simeq R\epsilon_*p_n^*(A_X),\;\hbox{and}\,\;
R\epsilon_*q_n^*(\bar A_n)\simeq \Hat\pi^*R\pi_{n*}(\bar A_n).$$
 The
object $R\mu_*(A)\in D_G^b({\rm pt})$ is defined by the sequence
of triples
$$\{(R\mu_*(A_X),\,R\pi_{n*}(\bar A_n),\,e_n)\},$$
 where $e_n$ is given by
the composition of the the above isomorphisms
 with
$$R\epsilon_*(a_n):R\epsilon_*p_n^*(A_X)\to R\epsilon_*q_n^*(\bar
A_n).$$ In this way, we have defined a functor
$R\mu_*:D^b_{\varphi}(X)\to D^b_G({\rm pt})$.

 From Proposition \ref{PropEnumerate}, one deduces  the following proposition.
\begin{Prop}\label{PropIsomrphismLambda2}
If $\vartheta$ is a $(\varphi,\varphi')$-cohomological morphism and
  \begin{enumerate}
 \item  $X$ is a  subanalytic space,
 \item  $G$ is a Lie group with a
finite number of connected components.
\end{enumerate}
 Then for any ${\mathcal A}\in {\mathfrak D}^b(X,\varphi)$, the objects
  $R\mu'_*(\delta'(\vartheta^{-1})^*({\mathcal A}))$ and
 $R\mu_*(\delta({\mathcal A}))$ of $D^b_G({\rm pt})$ are canonically isomorphic
 $$ R\mu'_*(\delta'(\vartheta^{-1})^*({\mathcal A}))   \simeq R\mu_*(\delta({\mathcal A})),$$
where ${R\mu_*}$ and ${R\mu'_*}$ are the functors induced by the
constant map $\mu:X\to{\rm pt}$.
 \end{Prop}

\subsection{Homotopic $G$-actions}
Let $\{ \varphi^s:G\to{\mathcal H}\,|\,s\in[0,\,1] \}$ be a
homotopy consisting of continuous group homomorphisms. We put
$\varphi:=\varphi^0,$  $\varphi':=\varphi^1$ and we write
$\varphi\sim_h\varphi'$.


 The maps $\phi_n,\phi_n':EG(n)\rightrightarrows E{\mathcal H(n)}$
are also homotopic by means of the family
\begin{equation}\label{ObviousHomotopy}
 \{\langle g,t\rangle\mapsto \langle \varphi^s(g),t\rangle
 \}_{s\in[0,\,1]}.
\end{equation}

In the same way, the
maps $\Phi_n,\Phi_n':BG(n)\rightrightarrows B{\mathcal H}(n)$  are
homotopic   through the homotopy $F_n$ defined by
 \begin{equation}\label{ObviousHomotopy1}
 BG(n)\times [0,\,1]\to B{\mathcal H}(n):\;\;
 ([g,t],\,s)\mapsto [\varphi^s(g),t]=:\Phi^s_n([g,t]),
  \end{equation}
 where $[g,\,t]$ is the image of  $\langle g,\,t\rangle$ by the projection to $BG(n)$.

  Since $\Phi_n$ and $\Phi'_n$ are homotopic, the ${\mathcal H}$-principal bundles
  $X(\varphi)(n)$ and   $X(\varphi')(n)$  are isomorphic.  In the following proposition, we recall  the construction of this isomorphism using the fold joint structure of $EG(n)$. Some stuff of this construction will be used for proving Proposition \ref{CompatibiXi}.

 \begin{Prop}\label{PropLambda0}
 Let $\varphi,\varphi':G\to {\mathcal H}$ be continuous group homomorphisms that are homotopic through a family
 of group homomorphisms. Then, for each
  positive integer $n$, there exist bundle isomorphisms over the identity of $BG(n)$
 $$\lambda_n: \Phi^*_n(E{\mathcal H}(n))\to \Phi'^*_n(E{\mathcal
H}(n))\;\;\hbox{and}\;\; \xi_n: X(\varphi)(n) \to X(\varphi')(n),$$
   which satisfy
   $\xi_n\circ Q_n=Q'_n\circ(\lambda_n\times 1).$
  \end{Prop}

{\it Proof.} Fixed the integer $n$,  we   omit in this proof the label $n$ in our notations; so we write $EG$ for $EG(n)$, $BG$ for $BG(n)$, $\Phi$ for $\Phi_n$, etc.

The fibre
${\mathcal H}$-bundles on $BG\times I$, $F^*(E{\mathcal H})$ and
$\check\pi^*(\Phi'^{*}(E{\mathcal H}))$ (where $\check \pi$ is the
projection of $BG\times I$ on the first factor) are isomorphic on
$BG\times\{1\}$.

 According to the proof of Theorem 9.6 of \cite{Hus} (see \cite[page 50]{Hus}),
 to extend the isomorphism to
 $BG\times I$, are sufficient  a covering $\{V_i\}_i$ of  $BG$,  trivializations of $F^*(E{\mathcal H})$ over the $V_i$'s
 and a family of maps $w_i:BG\to I$, with $w_i^{-1}(0,\,1]\subset V_i$.

\smallskip
{\it (i) Covering of $BG$.}
 For $i\in \{1,\dots,n\}$, let $t_i$ be the natural map $\langle g,t\rangle\in EG\mapsto t_i\in [0,\,1]$ (see  \cite[page 53]{Hus}).
 The family
 \begin{equation}\label{Coveringn}
 \{V_i:= \Hat\pi(t_i^{-1}(0,\,1 ])\}_{i=1,\dots,n},
  \end{equation}
  $\Hat\pi$ being the projection $EG\to BG$,    is an open covering of $BG$.

  \smallskip

{\it (ii) The functions $w_i$}. The maps $w_i:[g,\, t]\in BG\mapsto t_i\in[0,\,1]$ satisfy  the required property.

\smallskip

{\it (iii) Trivializations of $F^*(E{\mathcal H})$.}
Next, we construct a family $\{\eta_i\,|\, i=1,\dots,n\}$ of local sections
 of the ${\mathcal H}$-bundle $F^*(E{\mathcal H})$, where the  domain of  $\eta_i$ is $V_i\times I$.
 We put
 \begin{equation}\label{etai}
 \eta_i:V_i\times I\to F^*(E{\mathcal H}):\;\;
 ( [g,t],\,s)\mapsto\langle\varphi^s(g),
 t\rangle\varphi^s(g_i^{-1}).
  \end{equation}
More explicitly,
$$\langle\varphi^s(g),
 t\rangle\varphi^s(g_i^{-1}):=\big\langle
t_1\varphi^s(g_1)\varphi^s(g_i^{-1}),\dots,t_i1,\dots,t_{n}\varphi^s(g_{n})\varphi^s(g_i^{-1})\big\rangle.
$$
 The equality $[g,t]=[g',t']$ implies the existence of an element
 $b\in G$ such that $g'_j=g_jb$, for all $j$ such that $t_j\ne 0$.
 In particular, if
  $[g,t]=[g',t'] \in V_i$, then
   $b=g_i^{-1}g'_i$; thus, $\varphi^s(g_j)\varphi^s(g_i^{-1})=\varphi^s(g'_j)\varphi^s(g_i'^{-1})$, and
    the section $\eta_i$ is well defined.

\smallskip

With these ingredients, and following the construction given in \cite{Hus},
 one can define an isomorphism $\lambda$ of ${\mathcal H}$-bundles over the identity
 between $\Phi^*(E{\mathcal H})$ and $\Phi'^{*}(E{\mathcal H})$.

 The isomorphism $\lambda$ of principal bundles
 induces
another one between the associated bundles $X(\varphi)$ and
$ X(\varphi'),$
$$\xi:X(\varphi)\to
 X(\varphi'),$$ satisfying $\xi\circ Q=Q'\circ(\lambda\times 1).$

\qed

\noindent
 \subsection{ Compatibility of the $\xi_n$'s with each other.}\label{SubsectionCompatibility}
 The homotopy
 $$(g,s)\in G\times [0,\,1]\mapsto \varphi^s(g)\in{\mathcal H}$$
   between
 $\varphi$ and $\varphi'$ consisting of group homomorphims
  induces the homotopy $F_n$ between $\Phi_n$ and $\Phi'_n$, given (\ref{ObviousHomotopy1}). As
 the spaces $BG(n)$'s  and the $B{\mathcal H}(n)$'s are fold
 joints,
 the homotopies $F_n$ and $F_m$ are  compatible with the
 inclusions $BG(n)\hookrightarrow BG(m)$ and $B{\mathcal H}(n)\hookrightarrow B{\mathcal
 H}(m)$.
One has the following facts:
 \begin{enumerate}
 \item For $n<m$, the diagram
 $$\xymatrix{ F^*_n(E{\mathcal H}(n)) \ar@{^{(}->}[r] \ar[d] & F^*_m(E{\mathcal H}(m)) \ar[d]
 \\
 BG(n)\times [0,\,1]  \ar@{^{(}->}[r]  & BG(m)\times [0,\,1]
 }
 $$
 is commutative.
 \item The same property holds for the bundles
 $$\check\pi_n^*((\Phi_n')^*E{\mathcal H}(n))\;\; \hbox{and}\;\;
 \check\pi_m^*((\Phi_m')^*E{\mathcal H}(m)),$$
 where $\check\pi_n$ is the
 projection $BG(n)\times [0,\,1]\to BG(n).$
 \item We denote by   $\Hat\pi_{n+1}$   the quotient map $EG(n+1)\to BG(n+1)$, and for $j=1,\dots, n+1$ we put
   $V_j(n+1):=\Hat\pi_{n+1}(t_j^{-1}(0,\,1])$. Then, for $i=1,\dots, n$,  $V_i(n+1)\cap BG(n)=V_i$, where
 $V_i:=\Hat\pi_n(t_i^{-1}(0,\,1]),$ as in the proof of Proposition \ref{PropLambda0}.

\item  Analogously, for  $j=1,\dots,n+1$ one defines the maps
 $$w_{j}(n+1):BG(n+1)\to [0,\,1]:\;\, [g,\,t]\mapsto t_j.$$
 Obviously, if $i=1,\dots,n$ the restriction of $w_{i}(n+1)$   to $BG(n)$ coincides with the map $w_i$ defined in the proof of Proposition \ref{PropLambda0}.
 \item For $j=1,\dots,n+1$ we put
  \begin{gather}
 \eta_{j}(n+1):V_{j}(n+1)\times I\to F_{n+1}^*(E{\mathcal
 H}(n+1)),\notag\\
 ( [g,t],\,s)\mapsto\langle\varphi^s(g),
 t\rangle\varphi^s(g_{j}^{-1}). \notag
  \end{gather}
  For $1\leq i\leq n$, the restriction of $\eta_{i}(n+1)$ to $V_i\times I$ coincides with the map $\eta_i$ defined in (\ref{etai}).
  \end{enumerate}

 By these compatibilities, the homeomorphism
 $$\lambda_{n+1}:\Phi_{n+1}^*(E{\mathcal H}(n+1))\to \Phi'^*_{n+1}(E{\mathcal H}(n+1))$$
  (constructed following the process reminded in Proposition \ref{PropLambda0}) restricted to $\Phi_n^*(E{\mathcal H}(n))$ coincides with $\lambda_n$.

 Thus, the sequences of
homeomorphisms $\{\lambda_n\}_n$ and $\{\xi_n\}_n$ constructed in this way are consistent
with the inclusions. We have the following proposition.

\begin{Prop}\label{CompatibiXi}
 Let $\varphi,\varphi':G\to {\mathcal H}$ be continuous group homomorphisms that are homotopic through a family
 of group homomorphisms, then for each $n$ there exist
 homeomorphisms
 $$\lambda_n:\Phi_n^*(E{\mathcal H}(n))\to \Phi'^*_n(E{\mathcal
 H}(n)),\;\;\hbox{and}\;\;
 \xi_n:X(\varphi)(n)\to  X(\varphi')(n), $$
 satisfying
 $$
 i'_{mn}\circ\lambda_n=\lambda_m\circ i_{mn},\;\;\;\bar\psi'_{mn}\circ\xi_n=\xi_m\circ\bar\psi_{mn},$$
 for $n<m$.
\end{Prop}

From Proposition \ref{CompatibiXi} together with Proposition
\ref{PropLambda0}, it follows the corollary.

\begin{Cor}\label{CorCompati}
 The family  $\xi=\{\xi_n\}$ is a $(\varphi,\varphi')$-cohomological morphism.
\end{Cor}

 Thus, one has the corresponding
functor $\xi^*$, mentioned in Section \ref{S:intro}, that gives an
equivalence between the categories ${\mathfrak D}^b(X,\,
\varphi')$ and ${\mathfrak D}^b(X,\, \varphi)$.

\medskip

\noindent {\bf Proof of Theorem \ref{Thm:EquivCohom}.}
 The theorem is an immediate consequence of  Corollaries  \ref{CorCompati}
 and \ref{CorIguality}.

 \qed

As a corollary of Theorem \ref{Thm:EquivCohom},  we  deduce the
following result, which can also be proved directly without making
use of the derived categories.

\begin{Cor}\label{Thm:EquivR}
  If $\varphi$ and $\varphi'$ are
   group homomorphisms as in Theorem
  \ref{Thm:EquivCohom},  then
  $$H_{\varphi,c}(X;\,{\mathbb C})\simeq H_{\varphi',c}(X;\,{\mathbb
  C}).$$
\end{Cor}

 {\it Proof.}
 The result can be proved
applying Theorem \ref{Thm:EquivCohom}  to the object ${\mathcal
C}\in {\mathfrak D}^b(X,\varphi)$ defined in Subsection \ref{The category}.

\qed

\begin{Prop}
 Let $G$ be a topological group   and  $\varphi,\varphi':G\to
{\mathcal H}$ such that $\varphi\sim_h\varphi'$.
 If $X$   and $G$   are locally compact, and $X$ has finite cohomological dimension,
 then
$$H_{\varphi',c}(X; D_X)\simeq H_{\varphi,c}(X; D_X).$$
\end{Prop}

{\it Proof.} It is a consequence of Theorem \ref{Thm:EquivCohom} applied to the object ${\mathcal D}$ defined in (\ref{Ddualizing}).

\qed

\smallskip

\noindent
{\bf Proof of Theorem \ref{Thm:EquivCohom2}.}
 The equivariant cohomology of $X$, with respect to the $G$-action
$\varphi$, with coefficients
 in $A\in D^b_{\varphi}(X)$  is $H_{\varphi}(X;\,A):=H(R\mu_*A)$. Thus,
 Proposition \ref{PropIsomrphismLambda2} and Corollary \ref{CorCompati} imply  Theorem \ref{Thm:EquivCohom2}.

 \qed

From Theorem \ref{Thm:EquivCohom2}, we  deduce.
\begin{Cor}\label{Thm:EquivR2}
  Under the hypotheses of Theorem \ref{Thm:EquivCohom2}, the equivariant cohomologies $H_{\varphi}(X;\,{\mathbb C})$ and
   $H_{\varphi'}(X;\,{\mathbb C})$ are isomorphic.
\end{Cor}

{\it Proof.}  Theorem \ref{Thm:EquivCohom2} applied to the object
${\mathcal C}\in {\mathfrak D}^b(X,\varphi)$ gives the corollary.

\qed


\section{Localization}\label{SectLocalization}

 In this section, we assume that $G$ is
the {\em torus} $T=U(1)^n$ which acts on a $T$-{\em compactifiable
space} $X$. This hypothesis guarantees that there are finitely
many   different types of $T$-orbits, although $X$ is not
necessarily compact.

The $T$-equivariant cohomology with complex coefficients of a
point $H_T({\rm pt};\,{\mathbb C})=H(BT;\,{\mathbb C})$ can be identified to the
algebra ${\mathbb C}[{\mathfrak t}^*_{\mathbb C}]$,  of polynomials on the complexification ${\mathfrak
t}_{\mathbb C}$ of the Lie algebra of $T$. In this identification  the generators of the polynomial algebra
 are considered with even degree \cite[page 3]{A-B}.

Let $\Sigma$ denote the multiplicative subset of ${\mathbb
C}[{\mathfrak t}^*_{\mathbb C}]$ consisting of the non-zero
polynomials, and by $F$ we denote the fixed point set of the
$T$-action.
 For $A\in D^b_T(X)$, from the localization theorem
\cite[Sect. 6]{G-K-M} one deduces that the restriction map defines
an isomorphism
 $$H_T(X;\,A)_{\Sigma}\to H_T(F;\,A)_{\Sigma}$$
between the the corresponding localization modules.

As $T$ acts  trivially on $F$, $ET\times_T F=BT\times F$. Thus, for the particular
case of complex coefficients
\begin{equation}\label{F-free}
H_T(F;\,{\mathbb C})\simeq H(BT;\,{\mathbb
C})\otimes_{\mathbb C} H(F;\,{\mathbb C}).
 \end{equation}
 Hence,
 \begin{equation}\label{Sigma}
 H_T(X; \,{\mathbb C})_{\Sigma}\simeq{\mathbb C}(t^*_{\mathbb
 C})\otimes_{\mathbb C} H(F; \,{\mathbb C}),
  \end{equation}
 where ${\mathbb C}(t^*_{\mathbb C})$ is the field of rational
 functions on ${\mathfrak t}_{\mathbb C}$. In other words,
 $ H_T(X;\,{\mathbb C})_{\Sigma}$ is the result of the extension of scalars in $H(F;\,{\mathbb C})$ from ${\mathbb C}$
 to ${\mathbb C}(t^*_{\mathbb C})$.

\smallskip

\noindent
 {\bf Proof of Theorem \ref{Cohom=}.} The isomorphism  is a consequence of
Corollary  \ref{Thm:EquivR2}  together with (\ref{Sigma}) and the
fact that   the cohomologies
    are ${\mathbb C}$-vector spaces.

 \qed

\begin{Cor}\label{Corofinite}
Under the hypotheses of Theorem \ref{Cohom=},
 $$\bigoplus_{S\in{\mathfrak C}}H(S;\,{\mathbb C})\simeq
\bigoplus_{S'\in{\mathfrak C'}}H(S';\,{\mathbb C}),$$
 where ${\mathfrak C}$ and ${\mathfrak C'}$ are the
sets of connected components of $F$ and $F'$, respectively.
 In particular,
if $F$ and $F'$ are finite sets, then
$$\#F=\#F'.$$
\end{Cor}

\begin{Cor}\label{Hphi(F)=Hphi'(F')}
If the hypotheses of Theorem \ref{Cohom=} are satisfied, then
$$ H_{\varphi}(F;\,{\mathbb C})\simeq H_{\varphi'}(F';\,{\mathbb
C}).$$
\end{Cor}
{\it Proof.}
 From Theorem \ref{Cohom=}  together with (\ref{F-free}), it
 follows
$$H_{\varphi}(F;\,{\mathbb C})\simeq H(BT;\,{\mathbb
C})\otimes_{\mathbb C}H(F;\;{\mathbb C})
  \simeq H(BT;\,{\mathbb C})\otimes_{\mathbb C}H(F';\,{\mathbb C}) \simeq H_{\varphi'}(F';\,{\mathbb
  C}).$$

\qed

 Let $\varphi$ be a $T$-action on $X$ and $F$ the set of fixed points.
   The kernel of the restriction map
  $H_{\varphi}(X;\,{\mathbb C})\to H_{\varphi}(F;\,{\mathbb C})$ is a torsion $H(BT;\,{\mathbb C})$-submodule.
  We can consider the $H(BT;\,{\mathbb C})$-submodule $M$ of $H_{\varphi}(F;\,{\mathbb C})$ consisting
  of those elements which admit an extension to equivariant cohomology classes on $X$. We will show that under certain hypotheses
  $M$ is isomorphic to $H_{\varphi}(X;\,{\mathbb C})$.

   Since
  $BT$ is a simply connected space,  the $E_2$ page of the Leray-Serre cohomology spectral
   sequence of the fibration $X\to ET\times_T X\to BT$ is
  $$E_2^{p,q}=H^p(BT;\,{\mathbb C})\otimes_{\mathbb C}H^{q}(X;\,{\mathbb C}).$$

   If $H^{\rm odd}(X;\,{\mathbb C})=0$, as the cohomology of
  $BT$ vanishes in the odd degrees, the differential $d_r:E^{p,q}_r\to E^{p+r,q-r+1}_r$ is zero; so,
   the  spectral sequence
   collapses at the step
   $E_2$. In the terminology of \cite{G-K-M}, the $T$-space $X$ is equivariantly
formal. Thus,
 $$H_{\varphi}(X;\,{\mathbb C})=H(ET\times_T
X;\,{\mathbb C})=H(BT;\,{\mathbb C})\otimes_{\mathbb
C}H(X;\,{\mathbb C}).$$
 That is,  $H_{\varphi}(X;\,{\mathbb C})$ is a free finite generated ${\mathbb C}[{\mathfrak t}^*_{\mathbb C}]$-module.
  Therefore, $H_{\varphi}(X;\,{\mathbb C})$ can be identified
  with a submodule of $H_{\varphi}(F;\,{\mathbb C})$ (see \cite[page
  42]{G-K-M}); that is,  $M$ is isomorphic to $H_{\varphi}(X;\,{\mathbb C})$.
 Furthermore,  since the $T$-action on $F$ is trivial,
$H_{\varphi}(F;\,{\mathbb C})$ is
 a free ${\mathbb C}[{\mathfrak t}^*_{\mathbb C}]$-module. Thus,
   we have proved the proposition.
   \begin{Prop}\label{PropDivisor}
   If $H^{\rm odd}(X;\,{\mathbb C})=0$, then $M$ is a free
   ${\mathbb C}[{\mathfrak t}^*_{\mathbb C}]$-submodule of $H_{\varphi}(F;\,{\mathbb C})$
   isomorphic to $H_{\varphi}(X;\,{\mathbb C})$.
 \end{Prop}

 When $T=U(1)$,   ${\mathbb C}[{\mathfrak t}^*_{\mathbb C}]$ is a principal entire ring, and
 the free finite generated submodule $M$ of the  ${\mathbb
C}[{\mathfrak t}^*_{\mathbb C}]$-module
 $H_{\varphi}(F;\,{\mathbb C})$ has associated the corresponding elementary divisors.

\begin{Prop}\label{ElemenDiv}
  Let $\varphi$ and $\varphi'$ be homotopic $U(1)$-actions on a compactifiable space
  $X$, satisfying the hypotheses of Theorem \ref{Cohom=}.
  If $H^{\rm odd}(X;\,{\mathbb C})=0$, then the elementary divisors of $M$ and $M'$
  as submodules of $H_{\varphi}(F;\,{\mathbb C})$ and  $H_{\varphi'}(F';\,{\mathbb C})$ (resp.) coincide.
   \end{Prop}

{\it Proof.}
 From Corollary \ref{Thm:EquivR2}, Proposition \ref{PropDivisor} and
 Corollary \ref{Hphi(F)=Hphi'(F')}, we deduce the
 proposition.

   \qed


\end{document}